\documentclass[10pt,letterpaper]{article}
\usepackage[utf8]{inputenc}
\usepackage{amsmath}
\usepackage{amsfonts}
\usepackage{amssymb}
\usepackage{graphicx}
\usepackage{color}
\usepackage[normalem]{ulem}
\usepackage{algorithm,algorithmic,multirow}
\usepackage[left=2cm,right=2cm,top=2cm,bottom=2cm]{geometry}
\newenvironment{Proof}{\noindent{\sc Proof.}}{\qed}
\newtheorem{theorem}{Theorem}[section]
\newtheorem{lemma}{Lemma}[section]
\newtheorem{cor}{Corollary}[section]
\newtheorem{rem}{Remark}[section]
\newtheorem{definition}{Definition}[section]

\newcommand{\qed}{\hfill$\Box$\par\medskip}

\def\bhag#1{\noindent
\setcounter{equation}{0}
\section{#1}
}

\def\RR{{\mathbb R}}
\def\CC{{\mathbb C}}
\def\ZZ{{\mathbb Z}}

\def\PPI{{{\rm I}\kern-1pt\Pi}}

\def\a{\alpha}
\def\b #1;{{\bf #1}}
\def\x{{\bf x}}
\def\k{{\bf k}}
\def\y{{\bf y}}
\def\u{\mathbf{u}}
\def\w{{\bf w}}
\def\z{{\bf z}}

\def\j{\mathbf{j}}

\def\O{{\cal O}}

\def\C{{\mathcal C}}

\def\argmax{\mathop{\hbox{{\rm arg max}}}}

\def\be{\begin{equation}}
\def\ee{\end{equation}}
\def\bea{\begin{eqnarray}}
\def\eea{\end{eqnarray}}
\def\eref#1{(\ref{#1})}

\def\donchitre#1#2{\vskip 6.5cm\noindent
\parbox[t]{1in}{\special{eps:#1.eps x=6.5cm y=5.5cm}}
\hbox to 7cm{}\parbox[t]{0.0cm}{\special{eps:#2.eps x=6.5cm y=5.5cm}}}

\def\tn{|\!|\!|}

\def\bs#1{{\boldsymbol{#1}}}

\begin{document}
\title{A unified method for super-resolution recovery and real exponential-sum separation}
\author{Charles K. Chui\\
Department of Mathematics, Hong Kong Baptist University, Hong Kong 
\thanks{This author is also associated with the Statistics Department of Stanford University, CA 94305, and his research is partially supported by U.S. ARO Grant W911NF-15-1-0385 and a GRF grant from the Hong Kong Research Council, under Project HKBU\#12300917.}\\
\textsf{email:} ckchui@stanford.edu.  \\
H.~N.~Mhaskar\\
Institute of Mathematical Sciences, Claremont Graduate University,\\ Claremont, CA 91711.
\thanks {The research of this author is supported in part by ARO Grant W911NF-15-1-0385.}\\
\textsf{email:} hrushikesh.mhaskar@cgu.edu. }    
\date{}
\date{}

\maketitle

\begin{abstract}
In this paper, motivated by diffraction of traveling light waves, a simple mathematical model is proposed, both for the multivariate super-resolution problem and the problem of blind-source separation of real-valued exponential sums. This model facilitates the development of a unified theory and a unified solution of both problems in this paper. Our consideration of the super-resolution problem is aimed at applications to fluorescence microscopy and observational astronomy, and the motivation for our consideration of the second problem is the current need of extracting multivariate exponential features in magnetic resonance spectroscopy (MRS) for the neurologist and radiologist as well as for providing a mathematical tool for isotope separation in Nuclear Chemistry. The unified method introduced in this paper can be easily realized by processing only finitely many data, sampled at locations that are not necessarily prescribed in advance, with computational scheme consisting only of matrix - vector multiplication, peak finding, and clustering. 
\end{abstract}

\bhag{Introduction}\label{intsect}

If a light source is bright enough, there is no difficulty in locating its position, even though it is extremely tiny. However, it would become a big challenge to tell if this source consists of light emanated from two or more point-masses. The reason is that light diffracts as it moves through space. More precisely, when light with wavelength $\lambda$, emanated from an extremely tiny point-mass, travels in a medium with refractive index $\gamma$ and converges to a spot with half-angle $\theta$, the image pattern of the spot is a bright circular spot, surrounded by rings of concentric circles, called the “Airy pattern”, with radius of the bright circular spot equal to
\be\label{abbe}
d_a = \frac{\lambda}{2(NA)},
\ee
where $(NA) = \gamma \sin\theta$ was coined ``numerical aperture'' by Ernst Abbe in his 1873 paper, where \eref{abbe} was derived. The radius $d_a$ in \eref{abbe} is called ``Abbe’s resolution barrier'', since light with wavelength $\lambda$ nanometers (nm), emanated from two point-masses that are less than $d_a$ nm apart, cannot be resolved from the overlapping Airy patterns, meaning that the centers of the two bright circular spots cannot be identified. In practice, even for a perfect lens with circular aperture, Abbe's barrier is perhaps a little bit too small. Indeed, within a few years, the 1904 Physics Nobel laureate, John Strutt, commonly called Lord Rayleigh, derived in his 1879 paper \cite{rayleigh1879xxxi} the following slightly larger resolution barrier:
 \be\label{rayleigh}
d_r = 1.22 \frac{\lambda}{2(NA)},
\ee
called the Rayleigh criterion, by considering light source with wavelength $\lambda$, emanated from two point-masses with the center of the Airy image pattern from one point-mass lying on the smallest surrounding concentric circle of the Airy image pattern of the other. In view of the resolution barrier  \eref{rayleigh}, even if an optical microscope with the highest available quality of lens elements is perfectly aligned and has the highest numerical aperture, the resolution remains to be limited to half of the wavelength of visible light, namely: from approximately 200 nm for violet light, 250 nm for green light, and 350 nm for red light, when the (human) visible light spectrum from 400 nm to 700 nm is considered.

For the past 130 years or so, active research and development effort in the advancement of optical microscopy has been on-going to overcome the resolution barrier governed by the Rayleigh criterion \eref{rayleigh}; and being able to break the 200 nm resolution barrier has been considered achieving “super-resolution”. The first notable achievement is a laser scanning fluorescence microscope, called 4$\pi$, that can  achieve super-resolution of 150 nm and even to 100 nm for violet light. Proposed in a German patent \cite{cramerpatent1972} and its follow-up scholarly publication \cite{cremer1974considerations}, a key idea of the inventors, C. Cremer and T. Cremer, was to create a perfect hologram that carries the whole field information of the emission of a point source in all  4$\pi$ directions. This approach was later improved by S. Hell, E. Stelzer, C. Lindek, and C. Cremer in \cite{hell1994confocal}. However, it must be pointed out that all previous and current modes of super-resolution imaging are based on fluorescence!

Fluorescence microscopy has become the most important physical phenomenon in modern biology and medicine. The feasibility of creating living materials with targeted expressible fluorescent proteins that revolutionized fluorescence microscopy, was motivated by the success of functional imaging by using Ca++ probes, leading to the method of morphological and structural imaging.  As the pioneer of Calcium imaging, Roger Tsien (1952 - 2016) also pioneered the use of light and color to “peek and poke” at living cells to study how they work, by developing a “rainbow” of probes based on jellyfish green fluorescent protein (GFP) (see \cite{tsien1998green}, \cite{haggerty_wall_street_j}). For this work, Tsien was awarded, along with O. Shimomura and M. Chalfie, the 2008 Nobel Prize in Chemistry, for the development of GFP imaging.  More recently, for the success of developing innovative methods and systems of fluorescence microscopy to achieve super-resolution, the 2014 Nobel Prize in Chemistry was awarded jointly to Eric Betzig, Stefan W. Hell, and William E. Moerner (see, for example, their more recent work in \cite{betzig2006imaging}, \cite{balzarotti2016nanometer}, and \cite{von2015correcting}).

Observing that the image of a tiny light source has an Airy pattern as a result of blurring by light diffraction and that the light intensity of the bright circular spot of the Airy pattern is well approximated by the Gaussian with appropriate “radius” $d_g$ (which can be approximated as closely as desired by adjusting its standard deviation), we propose, in this paper, to represent the image capture of any given light source, emanated from a collection of point-masses, as a linear combination of such Gaussians. Of course the size of $d_g$ depends on the specific application. For the classical optic problem as described above, we may set $d_g = d_r$, and for application to resolving the super-resolution fluorescence microscopy problem, with $\ge 50$ nm resolution, say, we may select $d_g = \frac{d_r}{m}$, with $4\le m\le7$, where $d_r$ is the Rayleigh criterion in \eref{rayleigh}.  Then, the method and theory developed in this paper can be applied to locate the point-masses, by isolating the centers of the Gaussians and computing their locations. In other words, the innovative method developed in this paper can be applied to recover the imagery resolution. In particular, when applied to  fluorescence microscopy, the super-resolution can be resolved.  As to the application to observational astronomy, we may set 
$d_g = d_s$, where
\be\label{sparrow}
d_s = 0.94\frac{\lambda}{2(NA)},
\ee
is called the Sparrow criterion, derived by C.W. Sparrow \cite{sparrow1916spectroscopic} for the resolution barrier, defined as the minimum angular separation between two stars that can still be perceived as separate by an observer. The angular diameter of the Airy disk is determined by the aperture of the instrument. Sparrow's resolution limit is reached when the combined light from two overlapping and equally bright Airy disks is constant along a line between the central peak brightness of the two Airy disks.

Although the mathematical model (using Gaussian sums) proposed in this paper is intended to solving the super-resolution problem, particularly for fluorescence microscopy in the biology and medical fields, where only two and three dimensional image data are of interest, it also carries out to solving the blind-source separation problem of real exponential sums, with applications to data visualization and understanding of magnetic resonance spectroscopy (MRS) by the neurologist and radiologist, where image data in any dimension higher than $1$ are of interest. Since the Gaussian function is “blind” to dimensions, we will develop the theory, method, and computational algorithms for the proposed mathematical model for arbitrary dimensions. More precisely, we will consider the mathematical model:
\be\label{model}
G(\y, v)=\sum_{\ell=1}^L \frac{a_\ell}{(4\pi v^2)^{q/2}} \exp\left(-\frac{|\y-\y_\ell|^2}{4v^2}\right),
\ee
where $a_1, \cdots, a_L$ are real numbers and $\y$, $\y_1, \cdots, \y_L$ are in any desired $q$-dimensional space $\RR^q$, with $q\ge1$. 

When applied to the problem of resolution recovery, we consider a light source emanated from the linear combination: 
\be\label{targetmeasure}
\tau_I^*=\sum_{\ell=1}^L a_\ell\delta_{\y_\ell}.
\ee
of unknown point-masses $\y_1, \cdots, \y_L$, where $\delta_{\y_\ell}$ denotes, as usual, the Dirac delta distribution at the point $\y_\ell \in \RR^q$. Note that the model  function $G(\y,v)$ in \eref{model} is the integral convolution of $\tau_I^*$ in \eref{targetmeasure} with the normalized gaussian function:
\be\label{normalgauss}
g_v(\y)=\frac{1}{(4\pi v^2)^{q/2}}\exp\left(-\frac{|\y|^2}{4v^2}\right)
\ee
with standard deviation $\sqrt{2} v >0$, where the normalization is so chosen that the family $\{g(\y,v) :  v>0\}$ constitutes a positive approximate identity, so that the model function $G(\y,v)$ in \eref{model} converges to the light source $\tau_I^*$ in the distribution sense, as $v$ tends to zero. In this regard, for sufficiently small but fixed $v$, $G(\y,v)$ is a more realistic model of the light source than the distribution $\tau_I^*$ itself, since in all real-world applications, any point-mass in $\RR^q$, no matter how tiny, has positive measure. In addition, since diffraction of traveling light waves is a blurring process, the model function $G(\y,v)$ provides a very reasonable representation of the light intensity for the totality of the bright circular spots of the Airy patterns in the captured diffracted image. Note that in view of the increasing values of the coefficients in the Gaussian sum $G(\y,v)$ for smaller values of $v$, the bright circular spots of the Airy patterns are brighter when the blurring process is less severe. In this paper, we consider the more complete (super or standard) resolution recovery problem of finding the number $L$ of point-masses, computing their locations and corresponding light intensities, from finitely many data samples of  $G(\y,v)$, where $v$ is fixed, according to the desired resolution. In particular, the Super-resolution and Exponent Recovery algorithm (or Algorithm 1, to be stated in this section below) can be applied to resolve this more complete resolution recovery problem. In applying this algorithm to resolution recovery from the defected (blurred) imagery data, we point out that if the parameter $v$ of the model function $G(\y, v)$ exceeds the resolution barrier (collectively denoted by $d_g$, that stands for $d_r$,  $\frac{d_r}{m}$, and $d_s$, for standard resolution, super-resolution of $50$ nm, and observational astronomy, respectively, as discussed above), then the bright circular spots of neighboring Airy patterns cannot be resolved. 

For the problem of blind-source separation of the multivariate real-valued exponential sum:
\be\label{chuidata}
f(\y)=\sum_{\ell=1}^L b_\ell\exp(2\y_\ell\cdot \y), 
\ee
where $b_1, \cdots, b_L$ are real numbers and the unknown exponents $\y$, $\y_1, \cdots, \y_L$ can be treated as the point-masses of the measure $\tau_I^*$ in \eref{targetmeasure}. The reason is that the exponential sum $f(\y)$ in \eref{chuidata} can be changed to the model function $G(\y, v)$ for $v = 1/2$, namely:
\be\label{modelExS}
G(\y, 1/2) = \exp(- |\y|^2) f(\y), 
\ee
by setting the coefficients $a_\ell$ of $\tau_I^*$ in \eref{targetmeasure} to be
\be\label{notationchange}
a_\ell = \pi^{q/2}\exp(|\y_\ell|^2) b_\ell.
\ee 
Hence, again Algorithm 1 can be applied to $G(\y, v)$, by setting $v = 1/2$, to compute the exponents $\y_1, \cdots, \y_L$  and coefficients $a_1, \cdots, a_L$; from which the coefficients $b_\ell$ of $f(\y)$ can be recovered by applying \eref{notationchange}. That is, the blind-source exponential sum $f(\y)$ in \eref{chuidata} is separated.

In this paper, motivated by our work \cite{bspaper}, we propose a  method to recover the distribution $\tau_I^*$; that is, the number $L$ of terms, the coefficients $a_\ell$ and the points $\y_\ell$ in \eref{targetmeasure}, based on the action of a (Hermite)-weighted polynomial kernel on the scaled model function:  
\be\label{vriddence}
G(\x) = {(4\pi v^2)^{q/2}}G(\y/2v, v) = \sum_{\ell=1}^L a_\ell \exp(-|\x-\x_\ell|^2)
 \ee
with $\y_{\ell}$ in \eref{model} is replaced by $2v\x_{\ell}$, so that $v$ no longer exists in \eref{vriddence}. Similar to the algorithm developed in \cite{bspaper}, our approach is to define an operator $\mathcal{U}_n$ based on some parameter $n$, so that $\mathcal{U}_n(G)$ is equivalent to the action of an approximate identity on $\tau_I^*$. The distribution $\tau_I^*$ itself can then be recuperated, as in \cite{bspaper}, by using elementary techniques such as thresholding, peak finding, and clustering. Leaving the precise definitions of the various quantities involved to Section~\ref{mainsect}, the algorithm is summarized as Algorithm~\ref{algfigure} below.

\begin{algorithm}[h]\label{algfigure}
\caption{\textbf{SERA}: Super-resolution and Exponent Recovery}
\begin{algorithmic}[1]
\item[{\rm a)}] \textbf{Input:}  A subset $\C\subset \RR^q$ and the vector of values $\mathbf{g}=(\mathbb{G}(\y))_{\y\in\C}$ with $\mathbb{G}$ as in \eref{networkdef}.
\item[{\rm b)}]  \textbf{Tunable parameters:} $\lambda\in (0,1]$, $\mu>0$, and $\eta>0$. 
\item[{\rm c)}] \textbf{Pre-computation:} 
%
\STATE Solve the under-determined system of equations \eref{quadeqns} with $A=2/\sqrt{3}$. Determine $N$ so as to minimize either the residual error or the condition number of the system 
$$
\sum_{\y\in\C}w_\y\psi_\k(\sqrt{2A}\y)=(2A)^{-q/2}\delta_\k, \qquad |\k|_1\le 2N^2
$$
(cf. \eref{quadeqns}), or otherwise to minimize $\sum_{\y\in\C}|w_\y|$. Set $N=\lambda n$.
\STATE For $\x$ in a sufficiently dense subset $\mathcal{B}$ of $[-3\sqrt{3}N/2,3\sqrt{3}N/2]^q $, calculate the matrix 
$$
\mathbf{A}_{N;\x,\u}=w_\y\Phi_N^*(\x,\y), \quad \mathbf{A}_{n;\x,\y}=w_\y\Phi_n^*(\x,\y), \qquad \x\in\mathcal{B}, \ \y\in\C.
$$
\item[{\rm d)}] Compute 
$$
\mathcal{U}_N(\nu;\mathbb{G})(\x)=(\mathbf{A_Ng})_\x,\quad \mathcal{U}_n(\nu;\mathbb{G})(\x)=(\mathbf{A_ng})_\x, \qquad \x\in\mathcal{B}, \ \u\in\C.
$$
and find $\mathcal{G}_N=\{\x\in \mathcal{B} : |\mathcal{U}_n(\nu;\mathbb{G})(\x)|\ge A_2\mu/2\}$.
\item[{\rm e)}] Use clustering to partition $\mathcal{G}_N$ into sets $\mathcal{G}_\ell$,  satisfying (ii) and (iii) of Theorem~\ref{maintheo}(a). Let $L$ be the number  of these clusters.
\item[{\rm f)}] Compute $\widehat{\x_{N,n,\ell}}$ as in \eref{estcenter}, $\hat{a_\ell}=\mathcal{U}_n(\nu;\mathbb{G})(\widehat{\x_{N,n,\ell}})$, $\ell=1,\cdots,L$.
\item[{\rm g)}] \textbf{Output:} $L$, $\widehat{\x_{N,n,\ell}}$ (approximation to $\y_\ell$), $\widehat{a_\ell}$ (approximation to $a_\ell$). 
\end{algorithmic}
\end{algorithm}

\begin{rem}\label{normalizedserarem}
{\rm
If the input to Algorithm~\ref{algfigure} is of the form $G(\y,v)$ as in \eref{model} rather than $\mathbb{G}(\y)$ as in \eref{networkdef}, then we need to make the following modifications.
\begin{enumerate}
\item In item a), we need to input the known parameter $v$.
\item Before c)(1), we need to add the following step: Set $\x=\y/2v$, $\y\in\tilde{C}$, denote the set of $\x$'s by $\C$ and the vector of values $(\mathbb{G}(\x))_{\x\in\C}$ by $\mathbf{g}$.
\item Replace step g) by the following:
\begin{enumerate}
\item[{\rm g)(1:)}] Set $\widehat{\z_{N,n,\ell}}=v\widehat{\x_{N,n,\ell}}$, $\widehat{\alpha_\ell}=(\pi v^2)^{q/2}\hat{a_\ell}$.
\item[{\rm g)(2:)}] \textbf{Output:} $L$, $\widehat{\z_{N,n,\ell}}$ (approximation to $\y_\ell$), $\widehat{\alpha_\ell}$ (approximation to $a_\ell$). 
\end{enumerate} 
\end{enumerate}
\qed}
\end{rem}
\begin{rem}\label{SERAvalidity}
{\rm
The algorithm SERA described in Algorithm~\ref{algfigure} (with the modifications as in Remark~\ref{normalizedserarem}) applies to any fixed value of $v$ without any restrictions. However, for most applications, specific values of $v$ are of great interest. For example, in application to the decomposition of exponential sums $f(y)$ in \eref{chuidata}, we choose $v=1/2$ to obtain the result $G(\y, 1/2) = \exp(- |\y|^2) f(\y)$ in \eref{modelExS} for immediate application of our SERA. As to super-resolution recovery, since any point-mass, no matter how tiny, has positive measure in $\RR^q$, it is more realistic to use the mathematical representation 
\be\label{pointmass}
\chi_{v_0}(\y - \y_\ell) = (g_{v_0}*\delta_{\y_\ell})(\y)= g_{v_0}(\y - \y_\ell), 
\ee
for an appropriately small $v_0$, as opposed to the measure $\delta_{\y_\ell}$ itself, for the point-mass with size proportional $v_0$. Hence, because the size of the point-mass is often unknown, the optimal choice of $v$ in applying SERA may require some further work, but certainly must satisfy the condition $0< v \le d_g$ for resolving the light source 
\be\label{targetfunction}
F_{v_0}(\y)= \sum_{\ell=1}^L a_\ell\chi_{v_0}(\y - \y_\ell) 
\ee
instead of $\tau_I^*=\sum_{\ell=1}^L a_\ell\delta_{\y_\ell}$ in \eref{targetmeasure}, as discussed above. Here, as described above, $d_g$ stands for $d_r$,  $\frac{d_r}{m}$ or $d_s$. Application of super-resolution recovery to fluorescence microscopy and observational astronomy will be discussed in Sub-section~\ref{fluomicroscopy} and Sub-section~\ref{obsevastr}, respectively, of the final section.
\qed}
\end{rem}

\begin{rem}\label{heateqnrem}
{\rm
It is worthwhile to point out the well-known fact for the classical partial differential equation of isotropic heat diffusion with heat conductivity constant  $c>0$ (which states that at any $\y$ in $\RR^q$ and time instant $t > 0$, the partial derivative of the heat content $u(\y, t)$ with respect to the time variable $t$ is proportional to the Laplacian of $u(\y, t)$ with respect to the spatial variable $\y$, with $c$ as the constant of proportionality), that the solution, with initial heat content $\tau_I^*$  in \eref{targetmeasure} at $t=0$, is given explicitly by $$u(\y, t) =  G(\y, \sqrt{ct}),$$ with $v = \sqrt{ct}$ in \eref{model} for all $t > 0$. The interested reader is referred to  \cite[pp.~339-348] {chuijiang}.  While the heat diffusion is a forward problem in that the heat content as time $t = \frac{v^2}{c}$ increases, the resolution recovery problem is a (much harder) inverse problem in that the initial heat source $\tau_I^*$  is to be determined by information of the heat content $G(\y, \sqrt{ct}) = G(\y, v)$, for backward time $t = \frac{v^2}{c}$ travel (to near zero). In contrast to the traditional numerical PDE - based approaches to this problem, in which the time variable is reduced step by step, iteratively, to $0$, our algorithm gives a single shot (non-iterative) solution, directly from the observations of the solution at any fixed time instant $t$, or a fixed $v$.  
\qed}
\end{rem}

The organization of this paper is as follows. We will state our main results in Section~\ref{mainsect}. The quadrature formula needed for the discretization of the definition of $\mathcal{U}_n$, as an integral operator, is discussed in Section~\ref{quadsect}. The proofs are given in Section~\ref{pfsect}, after developing some background in Section~\ref{backsect}. Finally, applications to fluorescence microscopy, observational astronomy, magnetic resonance spectroscopy (MRS) and isotope separation will be briefly discussed in Section~\ref{applsect}.

\bhag{Main results}\label{mainsect}

Recall that in deriving the scaled model $G(\x)$ in \eref{vriddence}, we set $\x=\z/2v$, $\x_\ell=\y_\ell/2v$. With this change of point-mass notations, we find it convenient to introduce 
\be\label{targetmeasurebis}
\tau_I=\sum_{\ell=1}^L a_\ell \delta_{\x_\ell},
\ee
for the point-masses $\x_\ell$ to  replace $\tau_I^*$, introduced in \eref{targetmeasure} for the the point-masses $\y_{\ell}$. 

Since real-world data are usually not exact, it is necessary to introduce an error term. For this purpose, we will allow possible perturbation and consider the more general measure:
\be\label{startmeasure}
\tau=\tau_I +\tau_c,
\ee
where the additional measure $\tau_c$ is a complete, sigma finite, Borel measure with bounded total variation $\|\tau_c\|_{q,BV}$ on $\RR^q$. 
Thus, the problem which we wish to address in this paper is the following: Given information of the form
\be\label{networkdef}
\mathbb{G}(\x)=\sum_{\ell=1}^L a_\ell\exp(-|\x-\x_\ell|^2) + \mathcal{E}(\x),  \qquad \x\in\RR^q,
\ee
where
\be\label{epsdef}
\mathcal{E}(\x)=\int_{\RR^q}\exp(-|\x-\u|^2)d\tau_c(\u),
\ee
we wish to recuperate  $a_\ell$, $\x_\ell$ for $\ell=1,\cdots,L$; eqivalently, to recuperate the ideal target measure $\tau_I$ in \eref{targetmeasurebis}.

In Sub-section~\ref{contmainsect}, we introduce the Super-resolution and Exponent Recovery Operator (SERO), and formulate our first main theorem in the continuous setting; i.e., where the  $\mathbb{G}(\x)$ is assumed to be known for \textbf{all} $\x\in\RR^q$. In Sub-section~\ref{quadsect}, we discuss the case when $\mathbb{G}(\y)$ is given only for $\y\in \C$ for some finite set $\C\subset\RR^q$. This will involve the construction of certain quadrature formulas. 

\subsection{The Super-resolution and Exponent Recovery Operator}\label{contmainsect}
The first of our main theorems, Theorem~\ref{maintheo}, explains the recuperation of the target measure from the information $\mathbb{G}(\x)$, $\x\in\RR^q$. This will be done in a similar way as in our two previous work \cite{bspaper, hermite_recovery}. 
We need first to introduce some notation. 

Recall that for multi-integer $\k\in\ZZ^q_+$, the Hermite function $\psi_\k$ is defined via the generating function (cf. \cite[Section~5.5]{szego})
\be\label{hermitegenfunc}
\sum_{\k\in\ZZ^q_+}\frac{\psi_\k(\x)}{\sqrt{2^{|\k|_1} \k!}}\w^\k=\pi^{-q/4}\exp\left(-\frac{1}{2}|\x-\w|^2+|\w|^2/4\right), \qquad \w\in \CC^q.
\ee
In the sequel, we let  $H: [0,\infty)\to [0,1]$ to be a fixed, infinitely differentiable, non-increasing, even function,   satisfying $H(t)=1$ for $t\le 1/2$ and $H(t)=0$ for $t\ge 1$.   

Using $H$ and the functions $\psi_\j$, we define two kernels:
\bea\label{kerndef}
\Phi_n(\x,\y)&=&n^{-q}\sum_{\j\in\ZZ^q_+} H\left(\frac{\sqrt{|\j|_1}}{n}\right)\psi_\j(\x)\psi_\j(\y), \nonumber\\
 \Phi_n^*(\x,\y)&=&\left(\frac{2}{n^2\pi}\right)^{q/2}\sum_{\j\in\ZZ^q_+} H\left(\frac{\sqrt{|\j|_1}}{n}\right)3^{|\j|_1/2}\psi_\j(\x)\psi_\j(2\y/\sqrt{3}), \quad n>0, \x,\y\in\RR^q.
\eea

Although only an integral operator is discussed in this section, we prefer to define the Super-resolution and Exponent Recovery Operator (SERO) in a more general setting using an integral with respect to a more general measure. 
Throughout this paper, the term measure will mean a complete, sigma-finite, Borel measure on $\RR^q$, whether it is signed or positive. 
For a measure $\nu$, its total variation measure will be denoted by $|\nu|$. 
Of particular interest are the Lebesgue measure on $\RR^q$ and discretely supported measures of the form $\sum_{\y\in\C} w_\y\delta_\y$, where $\C\subset\RR^q$ is a finite set. 
The use of the Stieltjes integral $\int fd\nu$ to denote $\sum_{\y\in\C} w_\y f(\y)$ allows us to suppress cumbersome notation associated with the choice of the set $\C$ and the weights $w_\y$. More importantly, the measure will depend upon the choice of the parameter $n$. 
The measure notation  facilitates keeping track of the dependence of the various constants on this parameter and formulating the conditions on the choice of $\C$ as well as the weights in a concise manner without making tacit assumptions.

\begin{definition}\label{unopdef}
Let $\nu$ be a measure. Then the Super-resolution and Exponent Recovery Operator (SERO) is defined for $|\nu|$-integrable functions $f$ by
\be\label{intopdef}
\mathcal{U}_n(\nu;f)(\x)=\int_{\RR^q} f(\u)\Phi_n^*(\x,\u)\exp(-|\u|^2/3)d\nu(\u), \qquad n>0, \ \x\in\RR^q.
\ee
If $\nu$ is the Lebesgue measure on $\RR^q$, we will write $\mathcal{U}_n(f)$ in place of $\mathcal{U}_n(\nu;f)$.
\end{definition}

 As a consequence of Lemma~\ref{momentgenlemma} below, we will see that for $n>0$, $\x\in\RR^q$,
\be\label{untotaun}
\mathcal{U}_n(\mathbb{G})(\x)=\int_{\RR^q}\Phi_n(\x,\y)d\tau(\y)=\sum_{\ell=1}^L a_\ell \Phi_n(\x,\x_\ell) +\int_{\RR^q}\Phi_n(\x,\y)d\tau_c(\y)= \sum_{\ell=1}^L a_\ell \Phi_n(\x,\x_\ell) +\mathcal{U}_n(\mathcal{E})(\x), 
\ee
where $\mathcal{E}$ is defined as in \eref{epsdef}.
We will prove in Lemma~\ref{kernlemma} below that $\Phi_n(\x,\y)$ is an approximation to the Dirac delta $\delta_{\x-\y}$. As in \cite{hermite_recovery}, this will lead to an approximate recovery of the target measure from the function $\mathbb{G}$.

In the sequel, we use the following notation.
\be\label{Metcdef}
M=\sum_{_\ell=1}^L|a_k|,\ \mu=\min_{1\le _\ell \le L}|a_\ell|,\ \eta=\min_{1\le k\not=j\le L}|\x_\ell-\x_j|, \ B=\max_{1\le \ell\le L}|\x_\ell|_\infty.
\ee

Our first main theorem is the following, where the constants $A_2$, $\alpha$ are defined in Lemma~\ref{kernlemma}.

\begin{theorem}\label{maintheo}
Let $n\ge 1$, 
\be\label{largeset}
\mathcal{G}_n =\{\x\in \RR^q : |\mathcal{U}_n(\mathbb{G})(\x)|\ge A_2\mu/2\},
\ee
and $\gamma>0$ be as in \eref{gammadef}. Furthermore, assume that
\be\label{noisecond}
|\mathcal{U}_n(\mathcal{E})(\x)|\le A_2\mu/8, \qquad \x\in\RR^q.
\ee
Then for sufficiently large values of $n$, each of the following statements holds.
\begin{enumerate}
\item[{\rm (a)}]
 There exists a partition $\mathcal{G}_{n,\ell}$, $\ell=1,\ldots, L$,  of $\mathcal{G}_n$ such that
\begin{enumerate}
\item[(i)] For $\ell=1,\ldots, L$, $\x_\ell\in \mathcal{G}_{n,\ell}$.
\item[(ii)] For $\ell=1,\ldots, L$, $\mathsf{diam}( \mathcal{G}_{n,\ell})\le 2\gamma/n\le \eta/2$.
\item[(iii)] For $\ell, j=1,\ldots, L$, $\ell\not=j$, $\mathsf{dist} (\mathcal{G}_{n,\ell},\mathcal{G}_{n,j})\ge\eta/2$.
\end{enumerate}
\item[{\rm (b)}] Let 
\be\label{capNcond}
N\ge \max(1,\frac{2\gamma}{\a})n, 
\ee
and for $\ell=1,\cdots,L$,
\be\label{estcenter}
\widehat{\x}_{N,n,\ell}=\argmax_{\x\in \mathcal{G}_{N,\ell}}|\mathcal{U}_n(\mathbb{G})(\x)|.
\ee
Then
\be\label{centeresterror}
|\widehat{\x}_{N,n,\ell}- \x_\ell| \le 2\gamma/N \le \a/n,
\ee
and for a suitably small $\varepsilon$ (cf. \eref{epssmall})
\be\label{basicampfound}
|\mathcal{U}_n(\mathbb{G})(\widehat{\x}_{N,n,\ell})-a_\ell\Phi_n(\widehat{\x}_{N,n,\ell},\widehat{\x}_{N,n,\ell})| \le 5\varepsilon.
\ee 
\end{enumerate} 
\end{theorem}

We will show in Lemma~\ref{kernlemma} that there exists a positive constant $A_1$, depending only on $H$ and $q$ such that
$$
|\Phi_n(\x,\y)| \le A_1, \qquad \x,\y\in \RR^q.
$$
Therefore,  Lemma~\ref{momentgenlemma} leads to
\be\label{enest}
|\mathcal{U}_n(\mathcal{E})(\x)| =\left|\int_{\RR^q} \Phi_n(\x,\y)d\tau_c(\y)\right| \le A_1\|\tau_c\|_{q,BV}.
\ee

Thus, the condition \eref{noisecond} stipulates that the measure $\tau_c$ should not dominate the least of the terms $a_\ell\delta_{\x_\ell}$. If perturbation $\tau_c$ dominates some of the point measures, then of course, one does not expect a solution to the problem.

\subsection{Discretization}\label{quadsect}

We observe that evaluation of $\mathcal{U}_n$ involves evaluation of an integral. A straightforward discretization using Monte Carlo or quasi-Monte Carlo method will destroy the error estimates in Theorem~\ref{maintheo}. The most natural choice for the numerical integration is perhaps the Gaussian quadrature rule based on the zeros of Hermite polynomials. However, using this rule requires the knowledge of $\mathbb{G}$ at exactly these points. In practice, one cannot control the placement of points at which the values of $\mathbb{G}$ are available (the \textbf{sampling points}). Therefore, we need a quadrature formula, analogous to the Gaussian quadrature formula, but based on the sampling points.

As explained earlier, in order not to clutter our notation with unnecessary details of the quadrature formula in the statement of our result, we find it convenient to use a measure theoretic notation. To do so, we first need to introduce some further notation and conventions.

In the following, for $\lambda>0$, the symbol $\mathbb{P}_\lambda^q$ (respectively, $\mathbb{P}_\lambda^{q,\square}$) denotes the class of all polynomials in $q$ variables with total (respectively, coordinatewise) degree $<\lambda^2$. The symbol $\Pi_\lambda^q$ (respectively, $\Pi_\lambda^{q,\square}$) denotes the class of all functions of the form $\x\mapsto P(\x)\exp(-|\x|^2/2)$, $P\in \mathbb{P}_\lambda^q$ (respectively, $P\in \Pi_\lambda^{q,\square}$). 
It is not difficult to verify from the definition that for any $\k\in\ZZ_+^q$, the function  $\x\mapsto\psi_\k(\x)\exp(|\x|^2/2)$ is in  $\mathbb{P}_\lambda^q$ with $\lambda=\sqrt{|\k|_1}$, and in fact,
$$
\Pi_\lambda^q=\mathsf{span} \{\psi_\k: |\k|_1<\lambda^2\}, \quad P\in \Pi_\lambda^{q,\square}=\mathsf{span} \{\psi_\k: |\k|_\infty<\lambda^2\}
$$
We observe further that if $P, Q\in \mathbb{P}_\lambda^q$ (respectively, $P, Q\in \Pi_\lambda^{q,\square}$) then $PQ\in \mathbb{P}_{\sqrt{2}\lambda}^q$ (respectively, $PQ\in \Pi_{\sqrt{2}\lambda}^{q,\square}$).

We will use the following convention regarding constants.

\noindent\textbf{Constant convention:}\\
\textit{ The symbols $c, c_1,\cdots$ will denote generic positive constants depending only upon $H$, $q$, and the number $S$ to be introduced later. Their values may be different at different occurrences, even within a single formula.
 Constants denoted by capital letters will retain their values.
 }
 
Our description of the quadrature formula in the abstract is given in the following definition.

\begin{definition}\label{mzdef}
Let $A, n>0$. A complete, sigma finite, Borel measure $\nu$ on $\RR^q$ will be called a \textbf{Marcinkiewicz-Zygmund (MZ) quadrature measure of class $(q,A,n)$} (abbreviated $\nu\in \mathcal{M}(q,A,n)$) if there exist constants $c, c_1>0$  independent of $A$ or $n$, such that for every $P\in \mathbb{P}_{\sqrt{2}n}^q$,
\be\label{absquad}
\int_{\RR^q}P(\y)d\nu(\y)=\int_{\RR^q} P(\y)\exp(-A^2|\y|^2)d\y,
\ee
\be\label{absbv}
\|\nu\|_{q,BV} \le cn^q,
\ee
and 
\be\label{mzabs}
\int_{\RR^q}|P(\y)|d|\nu|(\y)\le c_1\int_{\RR^q} |P(\y)|\exp(-A^2|\y|^2)d\y.
\ee
\end{definition}

Thus, in Theorem~\ref{quadtheo} below, the measure that associates the weight $w_\y$ with each $\y\in \C$ is in $\mathcal{M}(q,A,n)$. 

We now discuss the discretized operator $\mathcal{U}_n(\nu;f)$, based on $\nu\in \mathcal{M}(q,2/\sqrt{3},\sqrt{2}n)$. As an application of Lemma~\ref{discmomentgenlemma}, we will see that if $\nu\in \mathcal{M}(q,2/\sqrt{3},\sqrt{2}n)$, then with $\mathcal{E}$ as in \eref{epsdef},
\be\label{discuntotaun}
\left|\mathcal{U}_n(\nu; \mathbb{G})(\x) - \sum_{\ell=1}^L a_\ell \Phi_n(\x,\x_\ell)-\mathcal{U}_n(\mathcal{E})(\x)\right| \le cn^{3q-2}3^{-n^2/2}\sum_{\ell=1}^L |a_\ell|.
\ee

The following theorem is the analogue of Theorem~\ref{maintheo} when the values of $\mathbb{G}$ are known only at the sampling points.

\begin{theorem}\label{maintheodisc}
Let $\nu\in\mathcal{M}(q,2/\sqrt{3},\sqrt{2} n)$. Then Theorem~\ref{maintheo} holds if $\mathcal{U}_n(\mathbb{G})$ is replaced by $\mathcal{U}_n(\nu;\mathbb{G})$ and the condition \eref{noisecond} is replaced by
\be\label{discnoisecond}
|\mathcal{U}_n(\mathcal{E})(\x)|\le A_2\mu/16, \qquad \x\in\RR^q.
\ee
\end{theorem}

\begin{rem}\label{movingmassremark}
{\rm
Let us point out an important consequence of this theorem for the determination of moving point masses and the corresponding coefficients. 
Observe that since the measure $\nu$ in Theorem~\ref{maintheodisc} is obtained independently of the target measures themselves, the same set of sampling points can be used to detect the target measures in an entire class of such measures. In particular, if the starting measures depend upon a time parameter, then the power spectra $|\mathcal{U}_n(\nu; \mathbb{G})(\x)|$ for all of these can be computed \textbf{in parallel}, and the extrema of these over the threshold as indicated can be obtained by curve extraction algorithms.
\qed}
\end{rem}

We now discuss the existence of MZ quadrature measures supported on an \textbf{arbitrary} set of points in $\RR^q$, subject to some conditions.
First, we need some preparation.

Since the kernel $\Phi_n^*$ is not symmetric, we find it convenient to formulate our quadrature formula in an apparently more general form than necessary.  For, $A, n>0$, $I_{A,n}=[-3n/A,3n/A]$. If $\C\subset I_{A,n}^q$ is a finite set, we define
 \be\label{meshnormdef}
\delta_{A,n}(\C)=\max_{\x\in I_{A,n}^q}\min_{\y\in\C}|\x-\y|_\infty, \qquad \eta_A(\C)=\min_{\y,\z\in\C, \ \y\not=\z}|\y-\z|_\infty,
 \ee  
and recall the invariance relation:
\be\label{meshinvrel}
A\delta_{A,n}(\C)=B\delta_{B,n}((B/A)\C), \qquad  A\eta_A(\C)=\eta_B((B/A)\C).
\ee
In order to apply Theorem~\ref{quadtheo}, we will require that the quantity $n\delta_{2/\sqrt{3}, n}(\C)$ for the set $\C$ of points at which $\mathbb{G}$ is observed to be sufficiently small. 

Theorem~\ref{quadtheo} below is motivated by the Gaussian quadrature formula based on the zeros of Hermite polynomials,  and assures the existence of quadrature formulas that are exact for integrating products of polynomials based on arbitrary sampling points. 

\begin{theorem}\label{quadtheo}
Let $n>0$ be sufficiently large, $A>0$, and $\C\subset I_{A,n}^q$ be a finite set. There exists $\beta>0$ (depending only on $q$) such that if $\delta_{A,n}(\C)\le \beta/(nA)$, then there exist weights $\{w_\y\}_{\y\in\C}$ that assure the following properties: 

For all $P\in \mathbb{P}_{\sqrt{2}n}^{q, \square}$, (in particular, for all products of polynomials in $\mathbb{P}_n^q$),
\be\label{quadrature}
\sum_{\y\in\C}w_\y P(\y)\exp(-A^2|\y|^2)=\int_{\RR^q} P(\y)\exp(-A^2|\y|^2)d\y,
\ee
\be\label{quadtotalvar}
\max_{\y\in \C}|w_\y| \le cA^q\delta_{A,n}(\C)^q, \qquad \sum_{\y\in\C}|w_\y| \le cn^q,
\ee
and
\be\label{mzineq}
\sum_{\y\in\C}|w_\y P(\y)|\exp(-A^2|\y|^2)\le cA^q\int_{\RR^q} |P(\y)|\exp(-A^2|\y|^2)d\y,
\ee
where all of the above constants, collectively denoted by $c$, are independent of $A$.
\end{theorem}

\begin{rem}\label{discunremark1}
{\rm
We note that the number of samples of $\mathbb{G}$ required to calculate $\mathcal{U}_n(\nu;\mathbb{G})$ in the above theorem is of order $\O(n^{2q})$.
\qed}
\end{rem}

\begin{rem}\label{discunremark2}
{\rm
A simple way to compute the quadrature weights $w_\y$ is to solve the under-determined system of equations
\be\label{quadeqns}
\sum_{\y\in\C}w_\y\psi_\k(\sqrt{2A}\y)=(2A)^{-q/2}\delta_\k, \qquad |\k|_1\le 2n^2.
\ee
Depending upon the application as well as computing resources, $n$ is determined so as to make either the residual error or the condition number of the system in \eref{quadeqns} small.  
For any set $\mathcal{B}\subset\RR^q$, we can precompute a matrix $\mathbf{A}$ indexed by $\mathcal{B}\times \C$, defined by
$$
\mathbf{A}_{\x,\y}=w_\y\Phi_n^*(\x,\y),  \qquad \x\in\mathcal{B}, \ \y\in\C.
$$
 The expression $\mathcal{U}_n(\nu;f)(\x)$, $\x\in\mathcal{B}$, can be implemented using a matrix vector multiplication $\mathbf{A}\mathbf{f}$ where $\mathbf{f}$ is the vector $(f(\y))_{\y\in\C}$.
\qed}
\end{rem}

\bhag{Background}\label{backsect}
For the convenience of the reader, we review some known facts regarding Hermite functions in this section.    

The univariate Hermite functions $\{\psi_j\}$ satisfy the Rodrigues' formula; i.e., \cite[Formula~(5.5.3)]{szego}
\be\label{hermitedef}
\psi_j(x)= \frac{(-1)^j}{\pi^{1/4}2^{j/2}\sqrt{j!}}\exp(x^2/2)\left(\frac{d}{dx}\right)^j (\exp(-x^2)), \qquad x\in\RR, \ j=0,1,\cdots.
\ee
We note the orthogonality relation for $j, k=0,1,\cdots$ (\cite[Formula~(5.5.1)]{szego}):
\be\label{uniortho}
\int_\RR \psi_j(x)\psi_k(x)dx=\left\{\begin{array}{ll}
1, &\mbox{ if $j=k$,}\\
0, &\mbox{otherwise}.
\end{array}\right.
\ee
We have the recurrence formula 
\bea\label{recurrence}
x\psi_{j-1}(x)&=&\sqrt{\frac{j}{2}}\psi_j(x) + \sqrt{\frac{j-1}{2}}\psi_{j-2}(x),\quad j=2,3,\cdots,\nonumber\\
&&\psi_0(x)=\pi^{-1/4},\ \psi_1(x)=\sqrt{2}\pi^{-1/4}x\exp(-x^2/2),
\eea
for $j=0,1,2,\cdots$ and $x\in\RR$ (cf. \cite[Formula~(5.5.8)]{szego}), and the Bernstein inequality (cf. \cite{freud1972direct}):
\be\label{unibern}
\|P'\|_{p,\RR} \le cn\|P\|_{p,\RR}, \qquad P\in \Pi_n^1, \ 1\le p\le \infty.
\ee
The infinite-finite range inequality \cite[Theorem~6.1.6, Theorem~6.2.4]{mhasbk} states that for any $\delta>0$, $1\le p\le\infty$, $P\in \Pi_n^1$,
\be\label{unirangeineq}
\|P\|_{p,\RR\setminus [-\sqrt{2}n(1+\delta), \sqrt{2}n(1+\delta)]} \le c_1e^{-cn^2}\|P\exp(-(\cdot)^2)\|_{p,[-\sqrt{2}n(1+\delta), \sqrt{2}n(1+\delta)]}.
\ee
The Mehler formula \cite[Formula~(6.1.13)]{andrews_askey_roy} states that
\be\label{mehler}
\sum_{j=0}^\infty \psi_j(y)\psi_j(z)r^j= \frac{1}{\sqrt{\pi (1-r^2)}}\exp\left(\frac{2yzr-(y^2+z^2)r^2}{1-r^2}\right)\exp(-(y^2+z^2)/2), \qquad y, z\in\RR, \ |r|<1.
\ee
It is proved in \cite{askey1965mean} that
\be\label{hermitebds}
|\psi_n(x)|\le\left\{\begin{array}{ll}
cn^{-1/4}, &\mbox{ if $|x|\le n^{1/2}(1-n^{-2/3})$,}\\
cn^{-1/12}, &\mbox{if $x\in\RR$.}
\end{array}\right.
\ee

The multivariate Hermite functions are tensor-product extensions of the univariate ones, namely:
\be\label{multihermitedef}
\psi_\k(\x)=\prod_{j=1}^q \psi_{k_j}(x_j).
\ee
We note that
\be\label{hermiteortho}
\int_{\RR^q} \psi_\j(\z)\psi_{\bs\ell}(\z)d\z =\delta_{\j,\bs\ell},\qquad \j,\bs\ell\in \ZZ^q_+.
\ee
The estimates \eref{hermitebds} imply that
\be\label{multihermitebd}
\max_{\x\in\RR^q}|\psi_\k(\x)| \le c, \qquad \k\in\ZZ_+^q.
\ee
The infinite-finite range inequality \eref{unirangeineq} takes the form:
for any $\delta>0$, $1\le p\le\infty$, $P\in\Pi_n^q$,
\be\label{rangeineq}
\|P\|_{p,\RR^q\setminus [-\sqrt{2}n(1+\delta), \sqrt{2}n(1+\delta)]^q} \le c_1e^{-cn^2}\|P\|_{p,[-\sqrt{2}n(1+\delta), \sqrt{2}n(1+\delta)]^q}
\ee
The multivariate version of the univariate Mehler formula \eref{mehler} takes the form
\be\label{multimehler}
\sum_{\j\in\ZZ^q_+} \psi_\j(\y)\psi_\j(\z)r^{|\j|_1}= \frac{1}{(\pi (1-r^2))^{q/2}}\exp\left(\frac{2\y\cdot\z r-(|\y|^2+|\z|^2)r^2}{1-r^2}\right)\exp(-(|\y|^2+|\z|^2)/2), \qquad \y, \z\in\RR^q, \ |r|<1.
\ee

The result from (\cite[Lemma~3.1]{hermite_recovery}) is a summary of some important properties of the kernel $\Phi_n$ in \eref{kerndef}.
 
\begin{lemma}\label{kernlemma} 
Let $S>q$ be an integer. There exist constants $A_1, A_2, C, C_1>0$ such that each of the following statements hold.
\begin{enumerate}
\item[{\rm (a)}] For $\x,\y\in\RR^q$, $n=1,2,\cdots$,
\be\label{hermite_localization}
|\Phi_n(\x,\y)| \le \frac{A_1}{\max(1,(n|\x-\y|)^S)}.
\ee
\item[{\rm (b)}] 
For  $n=1,2,\cdots$,
\be\label{philowbd}
|\Phi_n(\x,\x)|\ge A_2, \qquad |\x|_\infty\le Cn.
\ee
\item[{\rm (c)}] For $n\ge 1$, $|\x|_\infty, |\y|_\infty\le Cn$, we have
\be\label{summdiagstable}
|\Phi_n(\x,\x)-\Phi_n(\y,\y)|\le C_1n^{-q}|\x-\y|.
\ee
\item[{\rm (d)}]   There exists $\a>0$ such that
\be\label{schwarzcond}
0\le \Phi_n(\x,\y)\le \Phi_n(\y,\y), \qquad \x, \y\in\RR^q, \ |\x-\y|\le \a/n,\   |\y|_\infty\le Cn, \ n\ge 1.
\ee
\end{enumerate}
\end{lemma}

\begin{cor}\label{phinormcor}
For $\x\in\RR^q$, $n=1,2,\cdots$, $1\le p < \infty$
\be\label{phinintest}
\int_{\RR^q}|\Phi_n(\x,\y)|^pd\y \le cn^{-q}.
\ee
\end{cor}

\begin{Proof}\ 
Let $\x\in\RR^q$, and in this proof only, $B=\{\y\in\RR^q : |\x-\y|\le 1/n\}$. Then \eref{hermite_localization} shows that
\begin{eqnarray*}
\int_{\RR^q}|\Phi_n(\x,\y)|^pd\y &\le& A_1^p\left\{\int_B d\y +\frac{1}{n^{Sp}}\int_{\RR^q\setminus B}\frac{d\y}{|\x-\y|^{Sp}}\right\}\\
&\le& cA_1^p\left\{n^{-q} + \frac{1}{n^{Sp}}\int_{1/n}^\infty r^{q-Sp-1}dr\right\}=cA_1^pn^{-q}.
\end{eqnarray*} 
\end{Proof}

\bhag{Proofs}\label{pfsect}
In this section we prove both Theorem~\ref{maintheo} and Theorem~\ref{maintheodisc} by applying the fundamental theorem established in \cite[Theorem~2.1] {hermite_recovery}. For the convenience of the reader, we reproduce this theorem below.
\begin{theorem}\label{hermite_recovery_theo}
For $n\ge 1$, let
\be\label{generaltaun}
\mathbb{T}_n(\x)=\sum_{\ell=1}^L a_\ell\Phi_n(\x,\x_\ell)+E_n(\x),
\ee
and
\be\label{hermite_largeset}
\mathcal{G}_n =\{\x\in \RR^q : |\mathbb{T}_n(\x)|\ge A_2\mu/2\}.
\ee
We set
\be\label{gammadef}
\gamma=\max\left(1, \left(\frac{8A_1M}{A_2\mu}\right)^{1/S}\right),
\ee
where $A_1, A_2, S$ are as in Lemma~\ref{kernlemma}, and $M$, $\mu$ are as in \eref{Metcdef}.
We assume that
\be\label{hermite_noisecond}
|E_n(\x)|\le A_2\mu/8, \qquad \x\in\RR^q.
\ee
Then for sufficiently large values of $n$, each of the following statements holds.
\begin{enumerate}
\item[{\rm (a)}]
 There exists a partition $\mathcal{G}_{n,\ell}$, $\ell=1,\ldots, L$,  of $\mathcal{G}_n$ such that
\begin{enumerate}
\item[(i)] For $\ell=1,\ldots, L$, $\x_\ell\in \mathcal{G}_{n,\ell}$.
\item[(ii)] For $\ell=1,\ldots, L$, $\mathsf{diam}( \mathcal{G}_{n,\ell})\le 2\gamma/n\le \eta/2$, 
\item[(iii)] For $\ell, j=1,\ldots, L$, $\ell\not=j$, $\mathsf{dist} (\mathcal{G}_{n,\ell},\mathcal{G}_{n,j})\ge\eta/2$.
\end{enumerate}
\item[{\rm (b)}] Let 
\be\label{hermite_capNcond}
N\ge \max(1,\frac{2\gamma}{\a})n, 
\ee
and for $\ell=1,\cdots,L$,
\be\label{hermite_estcenter}
\widehat{\x}_{N,n,\ell}=\arg\max_{\x\in \mathcal{G}_{N,\ell}}|\mathbb{T}_n(\x)|.
\ee
Then
\be\label{hermite_centeresterror}
|\widehat{\x}_{N,n,\ell}- \x_\ell| \le 2\gamma/N \le \a/n,
\ee
and for a suitably small $\epsilon$ (cf. \eref{epssmall})
\be\label{hermite_basicampfound}
|\mathbb{T}_n(\widehat{\x}_{N,n,\ell})-a_\ell\Phi_n(\widehat{\x}_{N,n,\ell},\widehat{\x}_{N,n,\ell})| \le 5\epsilon.
\ee 
\end{enumerate} 
\end{theorem}
In Theorem~\ref{hermite_recovery_theo}, one needs $n$ to be large enough so that
\be\label{largelambda}
n\ge \max\left(1, 4\gamma/\eta, 2B/C, 4\gamma/\sqrt{C}, \left(\frac{A_2}{4C_1\gamma}\right)^{1/(q+1)}\right),
\ee
and
\be\label{epssmall}
\epsilon=\max(|E_n(\x_\ell)|, |E(\widehat{\x}_{N,n,\ell})|)+\frac{2^SA_1M }{(n\eta)^S}+\frac{2MC_1\gamma}{ n^{q+1}}\le \frac{\mu A_2}{4}.
\ee

The proof of Theorem~\ref{maintheo} in the present paper depends on the following lemma that connects the two kernels defined in \eref{kerndef}.

\begin{lemma}\label{momentgenlemma}
For $n>0$, $\x,\y\in\RR^q$,
\be\label{gausstophin}
\Phi_n(\x,\y)=\int_{\RR^q}\exp(-|\y-\u|^2)\Phi_n^*(\x,\u)\exp(-|\u|^2/3)d\u.
\ee
\end{lemma}

\begin{Proof}\ 
By completing squares, we deduce from \eref{multimehler}  that for $\y\in\RR^q$, $\k\in\ZZ^q_+$,
$$
\sum_{\k\in\ZZ^q}\psi_\k(\y)\psi_\k(\z)r^{|\k|_1}=\frac{1}{(\pi (1-r^2))^{q/2}}\exp\left(-\frac{1+r^2}{2(1-r^2)}\left|\y-\frac{2r}{1+r^2}\z\right|^2-\frac{1-r^2}{2(1+r^2)}|\z|^2\right).
$$
Taking $r=1/\sqrt{3}$ in this formula, we get
\be\label{mehlerspecial}
\sum_{\k\in\ZZ^q}\psi_\k(\y)\psi_\k(\z)3^{-|\k|_1/2}= \left(\frac{3}{2\pi}\right)^{q/2}\exp\left(-|\y-\frac{\sqrt{3}}{2}\z|^2\right)\exp(-|\z|^2/4).
\ee
In addition, from the orthogonality relation \eref{hermiteortho} and the substitution $\u=\sqrt{3}\z/2$, we deduce that
\be\label{hermite_recovery}
\psi_\k(\y)=3^{|\k|_1/2}\left(\frac{2}{\pi}\right)^{q/2}\int_{\RR^q} \exp(-|\y-\u|^2)\exp(-|\u|^2/3)\psi_\k(2\u/\sqrt{3})d\u.
\ee
The formula \eref{gausstophin} follows from here and the definitions \eref{kerndef}.
\end{Proof}

\begin{rem}\label{momentremark}
{\rm Let us point out a connection with our companion paper \cite{hermite_recovery}. Using \eref{hermite_recovery}, it is easy to see that
\bea\label{hermite_moments}
\hat{\tau}(\k) &=&\int_{\RR^q}\psi_k(\y)d\tau(\y)=3^{|\k|_1/2}\left(\frac{2}{\pi}\right)^{q/2}\int_{\RR^q} \left\{\int_{\RR^q}\exp(-|\y-\u|^2)d\tau(\y)\right\}\exp(-|\u|^2/3)\psi_\k(2\u/\sqrt{3})d\u\nonumber\\
&=&3^{|\k|_1/2}\left(\frac{2}{\pi}\right)^{q/2}\int_{\RR^q}\mathbb{G}(\u)\exp(-|\u|^2/3)\psi_\k(2\u/\sqrt{3})d\u.
\eea
Consequently, $\mathcal{U}_n(\mathbb{G})$ is the same as the same as the operator denoted in \cite{hermite_recovery} by $\mathcal{T}_n(\{\hat{\tau}(\k)\})$. \qed
}
\end{rem}
\noindent\textsc{Proof of Theorem~\ref{maintheo}.}\\
The definitions \eref{kerndef}, \eref{intopdef} and Lemma~\ref{momentgenlemma} lead to \eref{untotaun}. Therfore, $\mathcal{U}_n(\mathbb{G})(\x)$ has the form $\mathbb{T}_n(\x)$ as in Theorem~\ref{hermite_recovery_theo}, with $E_n(\x)=\mathcal{U}_n(\mathcal{E})(\x)$; and Theorem~\ref{maintheo} follows from Theorem~\ref{hermite_recovery_theo}, with $\epsilon$ given by \eref{epssmall} and $n$ large enough to satisfy \eref{largelambda} and \eref{epssmall}. \qed

The  following lemma is the analogue of Lemma~\ref{momentgenlemma}, which will be used in the proof of Theorem~\ref{maintheodisc} below. 

\begin{lemma}\label{discmomentgenlemma}
Let $n>0$, and $\nu\in \mathcal{M}(q,2/\sqrt{3},\sqrt{2}n)$. Then for $\x,\y\in\RR^q$,
\be\label{discgausstophin}
\left|\Phi_n(\x,\y)-\int_{\RR^q}\exp(-|\y-\u|^2)\Phi_n^*(\x,\u)\exp(-|\u|^2/3)d\nu(\u)\right| \le cn^{3q-2}3^{-n^2/2}.
\ee
\end{lemma}

\begin{Proof}\ 
In view of \eref{absquad}, and the fact that if $|\k|_1, |\j|_1 \le 2n^2$, then $\psi_k(2\y/\sqrt{3})\psi_\j(2\y/\sqrt{3})=P(\y)\exp(-4|\y|^2/3)$ for some $P\in \mathbb{P}_{2n}^q$, we obtain that
\be\label{pf1eqn1}
\int_{\RR^q}\psi_\k(2\u/\sqrt{3})\psi_\j(2\u/\sqrt{3})d\nu(\u)=
\int_{\RR^q}\psi_\k(2\u/\sqrt{3})\psi_\j(2\u/\sqrt{3})d\u
=(\sqrt{3}/2)^q\int_{\RR^q}\psi_\k(\z)\psi_\j(\z)d\z
=(\sqrt{3}/2)^q\delta_{\k-\j},
\ee
and hence, 
\bea\label{pf1eqn2}
\lefteqn{\int_{\RR^q}\left\{\sum_{\k\in\ZZ^q}\psi_\k(\y)\psi_\k(2\u/\sqrt{3})3^{-|\k|_1/2}\right\}\psi_\j(2\u/\sqrt{3})d\nu(\u)}\nonumber\\ &=&(\sqrt{3}/2)^q3^{-|\j|_1/2}\psi_\j(2\u/\sqrt{3}) +\int_{\RR^q}\left\{\sum_{\k: |\k|_1>2n^2}\psi_\k(\y)\psi_\k(2\u/\sqrt{3})3^{-|\k|_1/2}\right\}\psi_\j(2\u/\sqrt{3})d\nu(\u).
\eea
Using \eref{multihermitebd} and \eref{absbv}, we obtain for $ |\j|_1 \le 2n^2$
\begin{eqnarray*}
\lefteqn{\sum_{\k: |\k|_1>2n^2}|\psi_\k(\y)|\int_{\RR^q}|\psi_\k(2\u/\sqrt{3})3^{-|\k|_1/2}\psi_\j(2\u/\sqrt{3})|d|\nu|(\u)}\\
&\le& cn^q\sum_{\k: |\k|_1>2n^2}3^{-|\k|_1/2}
 \le cn^q\sum_{\ell= 2n^2}^\infty\ell^{q-1}3^{-\ell/2}\\
&\le& cn^{3q-2}3^{-n^2}.
\end{eqnarray*}
In view of \eref{mehlerspecial}, this implies that for $|j|_1<n^2$,
$$
\left|(2/\pi)^{q/2}\int_{\RR^q}\exp(-|\y-\u|^2)\exp(-|\u|^2/3)3^{|\j|_1/2}\psi_\j(2\u/\sqrt{3})|d\nu(\u)-\psi_\j(2\u/\sqrt{3})\right|\le cn^{3q-2}3^{-n^2/2}.
$$
Together with the definition \eref{kerndef} of $\Phi_n$, $\Phi_n^*$, this implies \eref{discgausstophin}.
\end{Proof}

\noindent\textsc{Proof of Theorem~\ref{maintheodisc}.}\\

In view of  Lemma~\ref{discmomentgenlemma}, we deduce that
$$
\left|\int_{\RR^q}\Phi_n(\x,\y)d\tau(\x)-\mathcal{U}_n(\nu;\mathbb{G})(\x)\right| \le cn^{3q-2}3^{-n^2/2}M.
$$
The equation \eref{untotaun} therefore leads to \eref{discuntotaun}; that is,
$$
\mathcal{U}_n(\nu;\mathbb{G})(\x)=\sum_{\ell=1}^La_\ell\Phi_n(\x,\x_\ell)+E_n(\x),
$$
where
$$
|E_n(\x)|\le |\mathcal{U}_n(\mathcal{E})(\x)|+cn^{3q-2}3^{-n^2/2}M.
$$
If \eref{discnoisecond} is satisfied, \eref{hermite_noisecond} is satisfied for sufficiently large $n$. In addition, with $\epsilon$ given by \eref{epssmall}, we choose $n$ large enough to satisfy \eref{largelambda} and \eref{epssmall}.
Thus, $\mathcal{U}_n(\nu;\mathbb{G})(\x)$ has the form $\mathbb{T}_n(\x)$ as in Theorem~\ref{hermite_recovery_theo}, with $E_n(\x)$ satisfying \eref{hermite_noisecond}. Therefore, Theorem~\ref{maintheo} follows from Theorem~\ref{hermite_recovery_theo}. \qed

In the proof of Theorem~\ref{quadtheo} to be given below, we will set $A=1/\sqrt{2}$, and the general theorem follows with a simple change of variables (cf. \eref{meshinvrel}). However, the proof requires the following preparation.  

Given the set $\C$ we may partition $I_{1/\sqrt{2},n}^q$ into congruent cubes with sides of length between $3\delta_{1/\sqrt{2},n}(\C)$ and $4\delta_{1/\sqrt{2},n}$. 
Clearly, each cube will contain some point of $\C$. For each cube, we pick the point closest to the center of the cube, and ignore the rest of the points in the cube. 
The resulting subset $\C'\subset \C$ has exactly one point in each cube, and $\eta(\C')\le 2\delta_{1/\sqrt{2},n}(\C')\le 3\delta_{1/\sqrt{2},n}(\C)\le 3\eta(\C')$. 
In Theorem~\ref{quadtheo}, we take $w_\y=0$ if $\y\in\C\setminus \C'$. 
In this way, we can rename $\C'$ as $\C$, the cube containing $\y\in\C$ as $J_\y$, and assume
\be\label{uniformity}
\eta(\C)\le 2\delta_{1/\sqrt{2},n}(\C)\le 4\eta(\C).
\ee

A critical step in the proof of Theorem~\ref{quadtheo} is the following lemma.

\begin{lemma}\label{prequadlemma}
With the set up as in Theorem~\ref{quadtheo}, let $\{J_\y\}_{\y\in\C}$ be a partition of $I_{1/\sqrt{2},n}^q$ that comprises congruent cubes with side $c\delta_{1/\sqrt{2},n}(\C)$, such that for each $\y\in\C$, $\{\y\}= J_\y\cap \C$  as described above. Then 
\be\label{strongmz}
\sum_{\y\in J_\y}\int_{J_\y}|P(\x)-P(\y)|d\x \le \frac{1}{8} \int_{\RR^q} |P(\x)|d\x, \qquad P\in \Pi_{\sqrt{2}n}^{q,\square},
\ee
and
\be\label{premzineq}
\frac{3}{4}\int_{\RR^q} |P(\x)|d\x \le  \sum_{\y\in \C}|J_\y||P(\y)| \le \frac{5}{4}\int_{\RR^q} |P(\x)|d\x, \qquad P\in \Pi_{\sqrt{2}n}^{q,\square}.
\ee
\end{lemma}

\begin{Proof}\ 
In the following,  we write $I^q=I_{1/\sqrt{2},n}^q$, $\delta=\delta_{1/\sqrt{2},n}(\C)$ and set 
\be\label{pf3eqn1}
V_n(x,y)=n\Phi_n(x,y), \qquad x, y\in \RR.
\ee
Then Corollary~\ref{phinormcor}  (with $q=1$)  and \eref{hermite_localization} imply that for $x, u\in\RR$,
\be\label{pf3eqn7}
|V_n(x,u)| \le cn, \qquad \int_\RR |V_n(x,u)|dx\le c.
\ee
For $\x=(x_1,\cdots,x_q)$ and $\u=(u_1,\cdots,u_q)$, we write (with an abuse of notation)
\be\label{pf3eqn2}
V_n(\x,\u)=\prod_{j=1}^q V_n(x_j,u_j).
\ee
It is easy to verify that $V_n\in \Pi_{2n}^{q,\square}$ in each of its arguments, and if $P\in \Pi_{\sqrt{2}n}^{q,\square}$, then
\be\label{pf3eqn3}
P(\x)=\int_{\RR^q} V_n(\x,\u)P(\u)d\u.
\ee

Let $P\in \Pi_{\sqrt{2}n}^{q,\square}$. Then \eref{pf3eqn3} implies that
\bea\label{pf3eqn4}
\sum_{\u\in \C} \int_{J_\y}|P(\x)-P(\y)|d\x &\le& \sum_{\u\in \C} \int_{J_\y}\left|\int_{\RR^q} V_n(\x,\u)P(\u)d\u-\int_{\RR^q} V_n(\y,\u)P(\u)d\u\right|\nonumber\\
&\le&  \|P\|_{1,\RR^q}\sup_{\u\in\RR^q}\sum_{\y\in\C} \int_{J_\y} |V_n(\x,\u)-V_n(\x,\u)|d\x.
\eea

Let $\u\in\RR^q$, and $Q_j=V_n(\cdot,u_j)$. Let $\x\in\C$ and $J_\y=\prod_{j=1}^q [a_j,b_j]$. For $\x\in J_\y$, we write
$$ 
\hat{\x}_j=(y_1,\cdots,y_{j-1}, x_j,\cdots,x_q), \quad j=2,\cdots,q,\ \hat{\x}_1=\x,\ \hat{\x}_{q+1}=\y. 
$$
Then it is easy to see using \eref{pf3eqn7} and the fact that the lengths of each side of $J_\y$ is $\sim \delta$, that
\bea\label{pf3eqn6}
\lefteqn{\int_{J_\y} |V_n(\x,\u)-V_n(\y,\u)|d\x\le \sum_{j=1}^q \int_{J_\y} |V_n(\hat{\x}_j,\u)-V_n(\hat{\x}_{j+1},\u)|d\x}\nonumber\\
& \le& 
 \sum_{j=1}^q \left\{\prod_{k=1}^{j-1} |Q_k(y_k)|(b_k-a_k)
\right\}\left\{(b_j-a_j)\int_{a_j}^{b_j}|Q_j'(u)|du\right\}\left\{\prod_{k=j+1}^q \int_{a_k}^{b_k} |Q_k(x_k)|dx_k\right\}\nonumber\\
&\le&
 c\delta\sum_{j=1}^q (n\delta)^{j-1}\left\{\int_{a_j}^{b_j}|Q_j'(u)|du\right\}\left\{\prod_{k=j+1}^q \int_{a_k}^{b_k} |Q_k(x_k)|dx_k\right\}.
\eea
Hence, using \eref{pf3eqn7} and the Bernstein inequality \eref{unibern} with $p=1$, we have 
$$
\sum_{\y\in\C}\int_{J_\y} |V_n(\x,\u)-V_n(\y,\u)|d\x\le c\delta\sum_{j=1}^q (n\delta)^{j-1}\int_\RR |Q_j'(u)|du \le c(n\delta)^q.
$$
Thus,  if $\delta\le c/n$ for a properly chosen $c$, then 
\be\label{pf3eqn5}
\sup_{\u\in\RR^q}\sum_{\y\in\tilde{\C}} \int_{J_\y} |V_n(\x,\u)-V_n(\x,\u)|d\x\le 1/8.
\ee
 In view of \eref{pf3eqn4}, we have proved \eref{strongmz}.

Next, by applying the infinite-finite range inequality \eref{rangeineq} (with $\sqrt{2}n$ in place of $\lambda$, $\sqrt{2}-1$ in place of $\delta$ there) we observe that
\bea\label{pf3eqn8}
\lefteqn{\left|\int_{\RR^q}|P(\x)|d\x-\sum_{\y\in\C}|J_\y||P(\y)|\right|
\le \int_{\RR^q\setminus I^q}|P(\x)|d\x + \left|\sum_{\y\in\C}\int_{J_\y}\left(|P(\x)|-|P(\y)|\right)d\x\right|}\nonumber\\
&\le& c_1\exp(-cn^2)\int_{\RR^q}|P(\x)|d\x +\sum_{\y\in\C}\int_{J_\y}|P(\x)-P(\y)|d\x.
\eea
Therefore, for sufficiently large $n$ and $\delta<c/n$ as required by \eref{strongmz}, we obtain
$$
\left|\int_{\RR^q}|P(\x)|d\x-\sum_{\y\in\C}|J_\y||P(\y)|\right|\le \frac{1}{4}\int_{\RR^q}|P(\x)|d\x.
$$
This implies \eref{premzineq}.
\end{Proof}

The proof of Theorem~\ref{quadtheo} is now standard.\\

\noindent
\textsc{Proof of Theorem~\ref{quadtheo}.}
In this proof, let $\C=\{\y_1,\cdots,\y_M\}$. We consider the space $U :=\Pi_{\sqrt{2}n}^{q,\square}$ and the sampling operator $\mathcal{S} : U \to \RR^M$ defined by
$\mathcal{S}(Q)=(Q(\y_1),\cdots,Q(\y_M))$; and define a norm 
$$\tn\z\tn=\sum_{j=1}^M |J_{\y_j}||z_j|, \qquad \z=(z_1,\cdots,z_M),$$
on $\RR^M$, and the functional
$$
x^*(\mathcal{S}(Q))=\int_{\RR^q} Q(\x)d\x.
$$
defined on the range of $\mathcal{S}$. The estimate \eref{premzineq} implies that 
\be\label{pf5eqn1}
|x^*(\mathcal{S}(Q))| \le \frac{4}{3}\tn \mathcal{S}(Q)\tn, \qquad Q\in U.
\ee
The Hahn-Banach theorem then yields an extension $X^*$ of $x^*$ to $\RR^M$ satisfying
\be\label{pf5eqn2}
|X^*(\z)| \le \frac{4}{3}\tn\z\tn, \qquad \z\in \RR^M.
\ee
We identify $X^*$ with $(W_1,\cdots,W_M)\in\RR^M$. Then the fact that $X^*$ is an extension of $x^*$ (and the definition of $x^*$) yields \eref{quadrature}. The estimate \eref{pf5eqn2} is equivalent to 
\be\label{pf5eqn3}
|W_k|\le |J_{\y_k}|, \qquad k=1,\cdots, M.
\ee
The first estimate in \eref{quadtotalvar} follows from the fact that $|J(\y_k)|\sim \delta_n(\C)^q$ for all $k$, and the second estimate follows from the fact that $M\sim n^{2q}$. The estimate \eref{mzineq} follows from \eref{premzineq} and \eref{pf5eqn3}. \qed

\vskip -0.5cm
\bhag{Applications}\label{applsect}

While the main objective of our previous paper \cite{hermite_recovery} is to recover the number $L$ of point-masses represented by $\x_{\ell}$, their positions in $\RR^q$, and the corresponding coeficients $a_{\ell}$ from the target measure $\tau_I$ defined in \eref{targetmeasurebis}, we have now developed a rigorous method along with the algorithm SERA in Algorithm~\ref{algfigure} by using the simple unified model function $G(\y, v)$ in \eref{model} and its scaled formulation $G(\x)$ in \eref{vriddence}, not only to achieve the same goal as \cite{hermite_recovery}, but alsos to determine the real exponents along with their corresponding coefficients $b_{\ell}$ from the blind-source function $f(\y)$ in \eref{chuidata}. An additional advantage of this unified model is that we may extend the target measure $\tau_I$ to the target function $F_{v_0}(x)$, with $\y$ and $\y_{\ell}$ replaced by $\x$ and $\x_{\ell}$ in \eref{targetfunction}, respectively. In particular, this facilitates the following discussions where point-masses have positive measures, even possibly with different values. 

\subsection{Fluorescence microscopy}\label{fluomicroscopy}
In \cite{hermite_recovery}, we discussed applications to counting red blood cells, as well as capturing shapes and colors of moving living cells, by considering suitable pairwise disjoint groups $X_k:=\{\x_{k, \ell} : \ell=1,\cdots, L_k\}$, where the $\x_{k, \ell}$ are selected among the $\x_{\ell}$, with corresponding coefficients $a_{k, \ell}$. Then if each group $X_k$ represents a different living cell, this approach allows a different $X_k$ with possibly different geometric shape and colors $a_{k, \ell}$ from the others. For example, if $X_{k_0}$ has a fractal-like shape and its nucleus (called the brain) is darker and more colorful than those of its neighbors $X_k$ (with $k$ different from $k_0$); and moreover, if the “size” of $X_{k_0}$, taken at incremental time instants, would increase, then $X_{k_0}$ most probably represents a cancer cell (see Remark~\ref{movingmassremark}). Unfortunately, since the size of each $X_k$ is zero, being just a set of points in $\RR^q$, its shapes and colors are unacceptable for visualization. In this paper, we replace each point $\x_{k, \ell}\in \RR^q$  by its corresponding true “pixel” $\chi_{v_0}(\x - \x_{k,\ell})$ introduced in \eref{pointmass}. Then depending on the minimum distance among the points that constitute $X_k$, a suitable choice of $v_0 > 0$ can be so chosen that the corresponding cell, represented by the union of the pixels $\chi_{v_0}(\x - \x_{k, \ell})$ for $\ell=1,\cdots, L_k$, with light intensity of each entire pixel  $\chi_{v_0}(\x - \x_{k,\ell})$ given by $a_{k, \ell}$, yields a significantly improved picture of the cell for visualization. In particular, if the coefficients $a_{k, \ell}$ are assigned various colors, the difference in darkness and colors of a cancer cell should clearly stand out from those of the normal ones. 

Recently, three-dimensional imaging with optical resolution as high as ~20 nm in the lateral direction and 40 - 50 nm in axial dimension has been achieved. The resolution of these super-resolution fluorescence microscopy techniques can reach the molecular scale, so that even the molecules within a cell can be separated by applying the algorithm SERA for the study of its structures and processes. Separated imagery features so obtained should enable scientists to directly visualize other biological samples at the nanometer scale and will complement the insights obtained through traditional molecular and cell biology approaches; thereby significantly expanding our understanding of molecular interactions and dynamic processes in living systems.

\subsection{Observational astronomy}\label{obsevastr}

The same approach of replacing points by pixels as discussed above applies to observing a galaxy and its stars. An advantage of this approach is that galaxies are identified by their shapes and colors. For shape classification, the most widely used scheme consists of spiral galaxies, elliptical galaxies, and irregular galaxies. Being the most common type, a spiral galaxy, such as our own Milky Way, is a rotating disk of stars and nebulae, surrounded by a shell of dark matter and with a bright central region at the core, called “galactic bulge”. On the other hand, the shape of an elliptic galaxy is ellipsoidal or ovoid, with size ranging from only a few thousand light-years to over hundreds of thousand light-years in diameter. Irregular galaxies have no particular shape. Full of gas and dust, most irregular galaxies are very bright. Those that are over 13 billion light-years away (implying that these galaxies are very young as we see them, and most probably have lots of star formation going on within them) are mainly irregular. As to color classification, old stars exist in “red” regions, since such stars have swollen and cooled, and emit “reddish” radiation. On the other hand, many young stars exist in “blue” regions, since such stars live fast and die young, consuming fuel at a high rate to maintain high temperatures that emit “blueish” hot radiation. 

In anticipation of the launch of the James Webb space telescope (called JWST or Webb) in October, 2018, with the capability of seeing stars that are 13.5 billion light-years away, and assuming that the Big Bang occurred 13.7 billion years ago,  it is exciting to have the opportunity to understand how stars that are almost 200 million light-years old are born. To identify such stars by their colors, note that when free protons capture free electrons in a cloud of ionized hydrogen, called an H-II (or H-two) region, light of various wavelengths, in “red/pink” color, is emitted, as electrons hop down through energy levels. H-II regions are ionized by ultraviolet radiation from hot stars, indicating the birth of new stars. As to the observation of change in size, we propose to apply our algorithm SERA, first to separate these young stars, and then to adjust the parameter $v_0$ for each of such stars, particularly of the very young ones, at incremental observation time instants to observe the rate and shape changes (see Remark~\ref{movingmassremark}). 

\subsection{Magnetic resonance spectroscopy}\label{mrs}

Magnetic resonance imaging (MRI) is based on nuclear magnetic imaging (NMR), which is a technique used by chemists and physicists to analyze and characterize small molecules in solid, liquid, and gel-like solutions. By performing a few additional procedures with MRI, particularly in shimming the magnetic field to correct its inhomogeneity by tuning the the x-y-z directions, magnetic resonance spectroscopy (MRS) imaging provides much more useful imaging information than MRI. For instance, for the radiologist, while MRI is used to identify anatomical locations of tumors, MRS can be used to determine tumor types and aggressiveness, as well as to distinguish between tumor recurrence and radiation necrosis. However, the price to pay is that shimming the magnetic field induces free inductive decay (FID) in that the receiver of the magnetic coil, with the rotating the components of the magnetization vector in the x-y plane, which crosses the coil loops perpendicularly. This causes the sinusoidal signal to decay exponentially with a time constant, called T2. Furthermore, back-to-back FID even causes the signal to increase and then to decrease again, both exponentially. Unfortunately, although MRS is a time - spatial domain operation, the common approach, particularly for in vivo measurement, is performed in the Fourier domain, which is not favorable to FID and back-to-back FID induced artifacts. By applying the algorithm SERA introduced in this paper locally before taking Fourier measurement, we believe that the FID artifacts can be eliminated, or at least mostly suppressed.

\subsection{Isotope separation}\label{isotopesep}

In basic chemistry, an atom consists of protons, electrons, and perhaps neutrons as well.  If there are neutrons in an atom, then the protons and neutrons cluster together in the central part of the atom, called the nucleus. While protons and electrons carry the same number of but opposite electric charges, neutrons carry no electric charge at all. An element is made up entirely of one type of atom, but its nucleus may or may not consist of neutrons. The periodic table is a tabular arrangement of elements, ordered by their atomic numbers; that is, the number of protons. So if the nucleus of an element consists of neutrons, this element is said to have different isotopes. For example, there are 3 isotopes of carbon, called carbon-12, carbon -13, and carbon-14, all with 6 protons but additionally 6, 7, and 8 neutrons in the nucleus, respectively, as well. The difference of C-14 from C-12 and C-13 is that C-14 is radioactive and in time, one of its neutron becomes a proton, by losing an electron, so that in some 5,700 years (called its half-life), half of the  C-14 become the stable nitrogen isotope, N -14. Other examples of isotopes include: uranium-238 that decays to lead-206, uranium-235 that decays to lead-207, potassium-40 that decays to argon-40, and rubidium that decays to strontium-87, but all with different half-lives. The elements with the most isotopes are cesium and xenon with 36 known isotopes. Applications of isotope separation include: separating uranium isotope to prepare enriched uranium for use as nuclear reactor fuel, separating hydrogen isotopes to prepare heavy water for use as moderator in nuclear reactors, and concentration of lithium-6 for use in thermonuclear weapons. The algorithm SERA, for the special case of dimension $q=1$ in \eref{chuidata}, should facilitate, at least in determining the number of radioactive isotopes to be separated from an element compound, by replacing the univariate exponent $2\y_{\ell}$ in \eref{chuidata}, with $ - \ln 2/ t_{\ell}$, where $t_{\ell}$ denotes the half-life of the radioactive $\ell^{th}$ isotope.


\end{document}